\documentclass[aap,MSNbibl,citesort,seceqn,dvips]{arximspdf}

\doi{10.1214/11-AAP823} %kopijuoti is PTS
\volume{22}
\issue{5}
\pubyear{2012}
\firstpage{2048}
\lastpage{2066}

\makeatletter
\newtheorem{thmm}{Theorem}[section]
\newtheorem{cor}[thmm]{Corollary}
\newtheorem{lem}[thmm]{Lemma}

\newproclaim{defn}{Definition}
\newproclaim{ass}{Assumption}
\def\lam{\lambda}

\def\Del{\Delta}

\def\al{\alpha}
\def\be{\beta}
\def\ga{\gamma}
\def\ep{\varepsilon}
\def\sg{\sigma}
\def\th{\theta}

\def\cC{\mathcal{C}}
\def\cH{\mathcal{H}}
\def\cP{\mathcal{P}}
\def\tE{E}
\def\td{\tilde{d}}
\def\cX{\mathcal{X}}
\def\mR{\mathbb{R}}
\def\mN{\mathbb{N}}

\def\tW{\tilde{W}}
\def\tG{\tilde{G}}
\def\tE{\tilde{E}}
\def\hr{\hat{r}}
\def\hW{\hat{W}}
\def\hG{\hat{G}}
\def\bB{\bar{B}}

\makeatother

\begin{document}
\begin{frontmatter}

\title{Nonuniform random geometric graphs with location-dependent radii}
\runtitle{Nonuniform random geometric graphs}

\begin{aug}
\author[A]{\fnms{Srikanth K.} \snm{Iyer}\corref{}\ead[label=e1]{skiyer@math.iisc.ernet.in}\thanksref{t1}}
\and
\author[B]{\fnms{Debleena} \snm{Thacker}\ead[label=e2]{thacker9r@isid.ac.in}}
\runauthor{S. K. Iyer and D. Thacker}
\thankstext{t1}{Supported in part by UGC SAP IV and grant from
DRDO-IISc program on Mathematical Engineering.}
\affiliation{Indian Institute of Science and Indian Statistical Institute}
\address[A]{Department of Mathematics\\
Indian Institute of Science\\
Bangalore, 560012\\
India\\
\printead{e1}} %adresu isvedimo komanda gale!
\address[B]{Theoretical Statistics\\
\quad and Mathematics Unit\\
Indian Statistical Institute 7\\
S.J.S. Sansanwal Marg\\
New Delhi, 110 016\\
India\\
\printead{e2}}
\end{aug}

% HISTORY:
\received{\smonth{1} \syear{2010}}
\revised{\smonth{7} \syear{2011}}

% ABSTRACT
%
\begin{abstract}
We propose a \textit{distribution-free} approach
to the study of random geometric graphs. The
distribution of vertices follows a Poisson
point process with intensity function $nf(\cdot)$, where
$n \in\mN$, and $f$ is a probability density function on $\mR^d$.
A vertex located at $x$ connects via directed edges to other vertices
that are within a \textit{cut-off} distance $r_n(x)$.
We prove strong law results for (i) the critical cut-off function
so that almost surely, the graph does not contain any node with
out-degree zero for sufficiently large $n$ and
(ii) the maximum and minimum vertex degrees.
We also provide a characterization of the cut-off function for which
the number of nodes with out-degree zero converges in distribution
to a Poisson random variable. We illustrate this result for a class
of densities with compact support that have at most polynomial rates
of decay to zero. Finally, we state a sufficient condition for an enhanced
version of the above graph to be almost surely connected eventually.
\end{abstract}

% KEYWORDS
%

\end{frontmatter}

%s1 ###
\section{Introduction and main results}\label{s1}
In this paper we study the asymptotic properties related to connectivity
of random geometric graphs where the underlying distribution of the
vertices may not be uniform. A random geometric
graph (RGG) consists of a set of vertices that are distributed in
space independently, according
to some common probability density function. The edge set of the
graph consists of the set of all pairs of points that are within a
specified cut-off distance. Our point of departure from usual random
geometric graphs is the specification of a cut-off function $r(\cdot)$
that determines the edge set. A directed edge exists, from a
vertex located at $x$ to another vertex located at $y$, provided
the distance between $x$ and $y$ is less than $r(x)$.

Our motivation for the study of such graphs comes from applications
in wireless networks. In models of wireless networks
as RGGs, the nodes are assumed to be communicating entities that are
distributed randomly in space according to some underlying density. Nodes
are assumed to communicate effectively with other nodes that
are within a cut-off distance, that is, proportional to the transmission
power. Hence, the transmission power has to be sufficiently large for
the network to be connected. However, nodes that are within each
other's transmission range interfere and thus cannot transmit simultaneously.
In order to maximize spatial reuse, that is, the
simultaneous use of the medium by several nodes to communicate,
the transmission power should be minimized.
Thus the asymptotic behavior of the critical radius of connectivity
in a random geometric graph as the number of vertices becomes large is
of considerable interest.

Often the nodes are assumed to be distributed in $[0,1]^d$
according to a Poisson point process of intensity $n \in\mN$.
In this case, it is known that the critical connectivity radius
scales as $O((\log n/n)^{1/d}).$ If the underlying density is nonuniform
but bounded away from zero, then (see Penrose~\cite{Pen99}) the asymptotic
behavior of the largest nearest neighbor distance in the graph is
$O((f_0^{-1}\log n/n)^{1/d}),$ where $f_0 > 0$ is the minimum of the
density over its support. Note that the asymptotics of the largest
nearest neighbor
distance is determined by the reciprocal of the minimum of the
density, since it is in the vicinity of the minimum that the nodes
are sparsely distributed.

In many applications such as
mobile ad-hoc networks and sensor networks,
the distribution of the nodes may be far from uniform (see, e.g., Foh
et al.~\cite{Foh}, Santi~\cite{San}). In the case
of nonuniform distribution of nodes, it is not efficient from
the point of view of maximizing spatial reuse, for all nodes to
use the same cut-off radius. Nodes near the mode of the density
require a much smaller radius than those at locations where
the density is small. Further, the infimum of the density over its support
could be zero ($f_0=0$). In such cases the asymptotics of the
largest nearest neighbor distance or the connectivity threshold
will be very different from that given above. One of the major
objectives in a wireless sensor network is to maximize battery life,
and hence it is important to minimize the energy expended in data
transmission. These considerations leads us to RGGs with
location-dependent choice of radii. In many applications,
it is assumed that the nodes know or can effectively estimate
their location
(Akkaya and Younis~\cite{Akk}, Langendoen and Reijers~\cite{Lang}).

We prescribe a formula for a critical location dependent cut-off
radius depending on the intensity, so that
almost surely the resulting graphs do not have isolated nodes eventually.
A useful property of the graphs we construct is that the distribution
of out-degree is independent of the location of the nodes. Vertex degree
distributions are important in designing algorithms for distributed
computations over wireless networks where the performance worsens with
increasing vertex degrees (Giridhar and Kumar~\cite{Giri}). We derive
strong\vadjust{\goodbreak} law bounds for the maximum and minimum vertex degrees.
By considering a finer parametrization of the cut-off function,
we show that the number of vertices with zero out-degree converges
to a Poisson distribution, under some conditions on the underlying
density. We illustrate the result with some examples. A result
of this nature for usual random geometric graphs with the uniform and
exponentially decaying densities of nodes can be found in Penrose~\cite{Pen97},
Gupta and Iyer~\cite{B}, respectively.

The solution to the connectivity problem for random geometric graphs
for nonvanishing densities with compact support and
dimensions $d \geq2$ can be found, for example, in Chapter 13
of Penrose~\cite{Pen03}. In one dimension, this problem is studied for
densities in $[0,1]$ with polynomial rate of decay to zero
in Han and Makowski~\cite{Mako}, while in Gupta, Iyer and
Manjunath~\cite{Iyer} the density is assumed to be exponential or
truncated exponential. In two dimensions, the asymptotic distribution
for the critical connectivity threshold for a large class of densities,
including
elliptically contoured distributions, distributions with independent
Weibull-like marginals and distributions with parallel level curves is
derived in Hsing and Rootzen~\cite{Hsing}. In dimensions $d \geq2$,
Penrose~\cite{Pen98} obtains the
asymptotic distribution for the connectivity threshold when the nodes
are distributed according to a standard normal distribution. In this
paper, we
derive a sufficient condition for the RGGs with location-dependent radius
to be almost surely connected eventually.

In summary, our primary motivation in proposing the study of graphs
with location
dependent radii is as follows. It is to enable the design and study of
wireless and sensor networks that allow nonstandard distribution
of nodes obtained by fitting densities to empirical data obtained from
actual deployments. Given such a density, each individual node can be programmed
to choose a transmission radius depending on its location so that the
network is connected with high
probability. Further, any change in the underlying distribution over
time due
to failures, re-deployments, etc., can be easily accommodated by appropriately
changing the transmission radii. As far as analyzing these graphs is concerned,
the key features to contend with are that the edges are directed and that
the cut-off radius is specified implicitly.

In order to state our results we need some notation.

%s1.1 ###
\subsection{Notation}
Let $f$ be a continuous\vspace*{-2pt} probability density function with
support $S \subset\mR^d$. For any random variable $X$ with density
$f$, we denote it by $X \stackrel{d}{\sim} f$.
Let the metric on $\mR^d$ be given by one of the $\ell_p$ norms
$1 \leq p \leq\infty$, denoted by $\Vert \cdot\Vert$. Let $\th_d$ denote
the volume of the unit ball in $\mR^d$. We denote by $\bar{B}$ the
closure of the set $B$.
Let $\cX= \{X_1,X_2, \ldots\}$ be a sequence of i.i.d. points
distributed according to $f$. Let $\{N_n\}_{n \geq1}$ be a nondecreasing
sequence of Poisson random variables with $E[N_n] = n$ and define
the sequence of sets
\begin{equation} \cP_n = \{X_1, X_2, \ldots, X_{N_n}\},\qquad n \geq
1. \label{eqnPn}
\end{equation}
Note that $\cP_n$ is a Poisson point process with intensity $nf$.
For any $r > 0$ and $x \in\mR^d$ we denote by
$B(x,r)$ the open ball of radius $r$ centered at $x$.\vadjust{\goodbreak}
For any Borel set $B \subset\mR^d$ and any point
process $\cP$, $\cP[B]$ represents the number of points of $\cP$ in~$B$,
and define
\begin{equation} F(B) := \int_B f(x) \,dx. \label{eqnF}
\end{equation}
We now define the random geometric graphs of interest with
location-depend\-ent radii.

\begin{defn}\label{defrgg} For any $n \geq1$, let $\cP_n$ be the set given by
(\ref{eqnPn}). For any function $r\dvtx \mR^d \to[0, \infty)$, the random
geometric graph $G_n(f,r)$ is defined to be the
graph with vertex set $\cP_n$, and the directed edge set
\[
E_n = \{ \langle X_i, X_j \rangle\dvtx X_i,X_j \in\cP_n, \Vert X_i -
X_j\Vert
\leq r(X_i) \}.
\]
\end{defn}

We will also consider an augmented version of the above random geometric
graph which is obtained by making all the edges in $G_n(f,r)$ bi-directional.

\begin{defn}\label{defrgge}
The enhanced random geometric graph $\tG_n(f,r)$
associated with the graph $G_n(f,r)$ is
defined to be the graph with vertex set $\cP_n$ and (undirected) edge
set
\[
\tE_n = \bigl\{ \{X_i,X_j\}\dvtx X_i,X_j \in\cP_n, \langle X_i, X_j \rangle
\in E_n \mbox{ or } \langle X_j, X_i \rangle\in E_n \bigr\}.
\]
\end{defn}

In the communication application described in the \hyperref[s1]{Introduction}, the
following procedure will give a graph whose edge set will \textit{contain} the
edges of the enhanced graph $\tG_n$. Upon deployment, the
nodes broadcast their radius. All nodes reset their transmission radius
to be the maximum of their original radius and the ones they receive from
the broadcast. Note that this is done only once. Clearly all the directed
links in the original graph now become bi-directional together with
possible creation of some directed edges. Thus if the enhanced graph
$\tG_n$ is connected, then so is the graph obtained by this procedure.
%In other words, for each node $X_i \in\cP_n$ define a
%new cut-off radius as
%
%
%s1.2 ###
\subsection{Main results}
For any fixed $c > 0$, define the sequence of \textit{cut-off} functions
$\{ r_n(c,x) \}_{n \geq1}$, via the equation
\begin{equation} \int_{B(x,r_n(c,x))} f(y) \,dy = c \frac{\log
n}{n},\qquad x \in S. \label{eqncutoff1}
\end{equation}
Later we will have occasion to take $c$ to be a function of $x$ and $n$
as well.
We will denote $G_n(f,r_n(c,\cdot))$ by~$G_n$ and the associated
enhanced graph by~$\tG_n$, when $c$, $f$ are fixed and $r_n$ is as defined in (\ref
{eqncutoff1}).
By the Palm theory for Poisson point processes (Theorem 1.6, \cite
{Pen03}), the expected
out-degree of any node in $G_n$ will be
\[
E[\operatorname{deg}(X_1)] = n \int_{B(x,r_n(c,x))} f(y) \,dy = c \log n,\vadjust{\goodbreak}
\]
which is the same as the vertex degree in the usual random geometric
graph defined on uniform points in the connectivity regime.
Let $P^x$ and $E^x$ denote the Palm distributions of $\cP_n$
conditional on a vertex located at $x$.
By the Palm theory for Poisson point processes,
the expected out-degree
of a node located at $x \in\mR^d$ in the graph $G_n$
will be
\[
E^x[\operatorname{deg}(x)] = c \log n.
\]
Thus the expected vertex degree of a node in $G_n$ does not
depend on the location of the node. In fact, the number of points
of $\cP_n \setminus\{x\}$ that fall in $B(x,r_n(c,x))$ under $P^x$
will follow a
Poisson distribution with mean $c \log n$.

%Let the functions $r_n = r_n(c, \cdot))$ be as defined in
%(\ref{eqncutoff1}).
Let $W_n=W_n(c)$ be the number of nodes in $G_n$ that have
zero out-degree, that is,
\begin{equation} W_n = \sum_{X_i \in\cP_n} 1_{\{\cP
_n[B(X_i,r_n(c,X_i))\setminus\{X_i\}] = 0\}}. \label{eqnWn}
\end{equation}
For each $n \geq1$, define
\begin{equation} d_n = \inf\{ c > 0 \dvtx W_n = 0 \}. \label{eqncriticald}
\end{equation}
In other words, $d_n$ is the critical cut-off parameter, that is, the
smallest $c$, so that
the graphs $G_n(f,r_n(c,\cdot))$ do not have any node with zero out-degree.
Our first result is a strong law for this critical cut-off parameter $d_n$.
\begin{thmm}\label{thmmdnstar} Let $d_n$ be the critical cut-off parameter as defined
in~(\ref{eqncriticald}). Then
almost surely,
\begin{equation} \lim_{n \to\infty} d_n = 1. \label{eqndnstar}
\end{equation}
\end{thmm}

Let $G_n = G_n(f,r_n(c, \cdot))$ be the random geometric graphs as
defined in Definition~\ref{defrgg} with $r_n$ as
in (\ref{eqncutoff1}). Consider the enhanced random geometric graph
$\tG_n = \tG_n(f, r_n(c, \cdot))$ associated with $G_n$; see
Definition~\ref{defrgge}.
%This graph is specified in Definition~\ref{defrgge}, where $\tr_n$
%is as
%defined in (\ref{eqntr}) with $r$ replaced by $r_n$.
We now state a strong law result for the critical cut-off parameter
to eleminate isolated nodes in the graph~$\tG_n$.
Let $\tW_n$ be the number of isolated nodes, that is, nodes
with degree zero, in the enhanced graph~$\tG_n$. Define
\begin{equation} \td_n := \inf\{c> 0\dvtx \tW_n = 0 \}. \label{eqntdn}
\end{equation}
Clearly $\td_n \leq d_n$ by construction.
The following theorem shows that the threshold required
to eleminate isolated nodes in the enhanced graph $\tG_n$ is the same as
for the graph $G_n$.
\begin{thmm}\label{thmmtdn}
Let $\td_n$ be as defined in
(\ref{eqntdn}). Then, almost surely,
\begin{equation} \lim_{n \to\infty} \td_n = 1. \label{eqnlimsuptdn}
\end{equation}
\end{thmm}

The exact asymptotics for the connectivity threshold for
random geometric graphs requires a lot of elaborate computations; see
Chapter 13,~\cite{Pen03}. We provide a sufficient condition that\vadjust{\goodbreak}
requires only a local computation at each node and makes use
of the connectivity threshold for the usual uniform random geometric graphs.

Let $X$ be a random variable with probability density function $f$
with support $S \subset\mR^d,$ $d \geq2$. Suppose that $f$ admits a mapping
$h\dvtx S \to\mR^d$ such that $h(X)$ is uniformly distributed on $[0,1]^d$.
For example, if the coordinates of $X$ are independently
distributed, then the coordinate mappings of $h$ will be the marginal
distributions. Recall
that $\th_d$ denotes the volume of the unit ball. For any
$\ep> 0$, let
$\{m_n\}_{n \geq1}$ be the sequence defined by
\[
m_n(\ep)^d = (1+\ep) \frac{m \log n}{n \th_d}, \qquad n \geq1,
\]
where
\begin{equation} m = \max_{0 \leq j \leq d-1} \frac{2^j (d-j)}{d}.
\label{eqnm}
\end{equation}
For any set $B$, we denote by $h(B)$ the image of the set $B$ under
$h$. Suppose that the
functions $r_n(c,x)$ are as defined in (\ref{eqncutoff1}).
Define the sequence of functions,
\begin{equation}\qquad c_n(\ep,x) := \inf\{ c\dvtx h(B(x,r_n(c,x))) \supset
B(h(x), m_n(\ep))\},\qquad x \in S, \label{eqncn}
\end{equation}
and $\{r_n(c_n, \cdot)\}_{n \geq1}$ to be functions on $S$ that
satisfy the equation
\begin{equation} \int_{B(x,r_n(c_n,x))} f(y) \,dy = c_n(\ep,x) \frac
{\log n}{n},\qquad n \geq1. \label{eqncutoff1n}
\end{equation}
Let $G_n(f,r_n(c_n))$ be the graphs defined
as in Definition~\ref{defrgg} with $r$ replaced by $r_n(c_n)$.
Let $\tG_n(f,r_n(c_n))$ be the enhanced graphs associated with
$G_n(f,\break r_n(c_n))$, that is, the graphs obtained by making the edges
in $G_n(f,r_n(c_n))$ bi-directional.
\begin{thmm}\label{thmmconnectivity} Let $X \stackrel{d}{\sim} f$, and suppose that $h(X)$ is
uniform on $[0,1]^d$, \mbox{$d \geq2$}. Then for any $\ep> 0$,
%let $c_n = c_n(\ep, \cdot)$, $n \geq1$, be as defined in (
%Then
almost surely, the sequence of enhanced
random geometric graphs $\tG_n(f, r_n(c_n))$ is connected for all
sufficiently large $n$.
\end{thmm}

In particular, the above result implies that
\[
P(\tG_n(f, r_n(c_n)) \mbox{ is connected }) \to1 \qquad\mbox{as } n
\to\infty.
\]
Our next result is on strong law asymptotics for the maximum and
minimum vertex degrees for the sequence of graphs $G_n$.
Let $H\dvtx[0,\infty)\longrightarrow[0,\infty)$ be defined by $H(0)=1$ and
%
%e1 ###
%
\begin{equation}
H(a)=1-a+a \log a,\qquad a > 0.
\end{equation}
The function $H$ has a unique turning point at the minima $a=1$. Let
$H_{+}^{-1}\dvtx\break [0,\infty)\rightarrow[1,\infty)$ be the inverse of
$H$ restricted to $[1,\infty)$ and $H_{-}^{-1}\dvtx[0,1]\rightarrow[0,1]$
be the
inverse of the restriction of $H$ to $[0,1]$.

\begin{thmm}\label{thmmmaxmindeg} For any $c > 0$, let $\Del_n = \Del_n(c)$ be the maximum
and $\delta_n = \delta_n(c)$ be the minimum vertex out-degree of
the graph
$G_n = G_n(f, r_n(c, \cdot)).$ Then with probability 1,
\begin{equation} \limsup_{n\to\infty} \frac{\Delta_{n}}{\log
n}\leq c H_{+}^{-1}( c^{-1} ). \label{eqnDeln}
\end{equation}
If $c < 1$, then $\delta_n \to0$ almost surely. If $c > 1$, then
with probability 1,
\begin{equation} \liminf_{n \to\infty} \frac{\delta_n}{\log n}
\geq c H_{-}^{-1}( c^{-1} ). \label{eqndeln}
\end{equation}
\end{thmm}

Note that (\ref{eqndeln}) does not shed any light on what happens to the
minimum vertex degree when $c=1$. This requires a finer parametrization
of the cut-off function. For any $\be\in\mR$, define
$\hr_n(\be,\cdot)\dvtx S \to[0,\infty)$ to be functions satisfying
\begin{equation} F(B(x,\hr_n(\be,x))) = \int_{B(x,\hr_n(\be,x))}
f(y) \,dy = \frac{\log n + \be}{n}, \label{eqnrPoi}
\end{equation}
for all n sufficiently large for which $\log n + \be> 0$, and arbitrarily
otherwise. Since~$\be$ is fixed, we will write $\hr_n(x)$ for $\hr
_n(\be,x)$.
For each $n \geq1$, define the sets
%
%e3 ###
%e2 ###
%
\begin{eqnarray}
A_n(x) & := & \{ y \in\mR^d \dvtx \Vert x-y\Vert \leq\hr_n(x) + \hr
_n(y) \},
\label{eqnAn} \\
\hat{A}_n(x) & := & \{ y \in\mR^d \dvtx \max\{\hr_n(x), \hr_n(y)\}
\leq\Vert x-y\Vert \leq\hr_n(x) + \hr_n(y) \},
\label{eqnhAn}
\\
\qquad K_n(x,y) &:=& B(y,\hr_n(y)) \setminus B(x,\hr_n(x)),\qquad
x,y \in S. \label{eqnKn}
\end{eqnarray}

\begin{thmm}\label{thmmPoicgs}
Let $\hW_n$ be the number of nodes of out-degree zero in
the graph $\hG_n = G_n(f, \hr_n(\be,\cdot)),$ where $\hr_n(\be
,\cdot)$
is as
defined in (\ref{eqnrPoi}).
%Let $S_1, \ldots,S_m$, be a finite partition of $S$.
Suppose that $f$ satisfies the following
two conditions for all $n$ sufficiently large:
\begin{longlist}[(1)]
\item[(1)] There exists a constant $\al\in(0,1)$, such that for all
$n$ sufficiently large
\begin{equation} \inf_{x \in S} \inf_{y \in\hat{A}_n(x)}
F(K_n(x,y)) \geq\al\biggl( \frac{\log n + \be}{n}\biggr),
\label{eqnc1}
\end{equation}
\item[(2)]
\begin{equation} \sup_{x \in S} F(A_n(x)) = o(n^{\al- 1})
\qquad\mbox{as } n \to\infty. \label{eqnc2}
\end{equation}
\end{longlist}
Then,
\begin{equation} \hW_n \stackrel{d}{\to} \operatorname{Po}(e^{-\be}) \label{eqnpoicgs}
\end{equation}
as $n \to\infty$, where $\operatorname{Po}(\lam)$ denotes a Poisson random variable
with mean $\lam.$
\end{thmm}

As noted earlier, Poisson approximation results for the number of
isolated nodes are available only for a\vadjust{\goodbreak} small class of distributions.
Condition (\ref{eqnc2}) can be replaced by the following
sufficient condition which, as will be shown below,
is easier to verify for some classes of densities.
\begin{equation} \sup_{x \in S} F(B(x, 2\hr_n(x))) = o(n^{\al- 1})
\qquad\mbox{as } n \to\infty. \label{eqnc3}
\end{equation}

%for all sufficiently large $n$.
%
\begin{thmm}\label{thmmPoi2} Suppose that the hypothesis (\ref{eqnc1}) of
Theorem \ref{thmmPoicgs}
and (\ref{eqnc3}) are satisfied, then (\ref{eqnpoicgs}) holds true.
\end{thmm}

Our final result is an illustration of the use of the above theorem to
a~class of densities with compact support that have at most polynomial rates
of decay to zero. Let $f$ be a continuous density with compact support
$S \subset\mR^d$, $d \geq2$. Let $B_i = B(x_i,r_i)$, $i=1,2, \ldots,k,$
be nonintersecting balls such that $f(x)=0$
on the boundary of these balls. For each $i=1,2, \ldots,k$, there exists
integers $m_i$ and $p_{ij}$, $j=1,2 , \ldots m_i$ and
constants $0 \leq\eta_i < r_i < \delta_i$ such that, either
\begin{eqnarray*}
f(y) &=& \sum_{j=1}^{m_i} A_{ij} (\Vert y-x_i\Vert - r_i)^{p_{ij}},\qquad
y \in B(x_i,\delta_i) \setminus B_i\quad \mbox{and}\\
f(y)& =& 0, \qquad y \in B_i \setminus B(x_i,\eta_i),
\end{eqnarray*}
or
\begin{eqnarray*}
f(y) &=& \sum_{j=1}^{m_i} A_{ij} (r_i - \Vert y-x_i\Vert)^{p_{ij}},\qquad
y \in B_i \setminus B(x_i,\eta_i) \quad\mbox{and}\\
f(y) &=& 0, \qquad y \in B(x_i,\delta_i) \setminus B_i,
\end{eqnarray*}
and $f(y) > 0$ elsewhere in $S$.
The density can vanish, for example, over water bodies which may
contain islands (which are being approximated by balls), and the density
decays to zero polynomially near the boundary of these balls. In
particular, $S$ itself could be taken to be a ball with the density
decaying polynomially to zero at the boundary of $S$ in a radially
symmetric fashion. We will denote by $\cH$ the set of all densities of
the above form. The class $\cH$ contains many of the standard distributions
such as the uniform, triangular, beta, etc., truncated versions of standard
distributions with unbounded support such as the Gaussian, gamma, etc. In
Santi~\cite{San} the density is assumed to be continuous and bounded
away from 0
over $[0,1]^2$, which is contained in $\cH$. The class $\cH$ also contains
the higher dimensional extensions of the polynomial densities considered
in Foh et al.~\cite{Foh} and Han and Makowski~\cite{Mako}. As
remarked earlier, since our motivation was to allow for nonstandard
densities, $\cH$ contains functions that are the modulus of analytic
functions in the interior of $S$.
\begin{cor}\label{corPoi} The conclusion of Theorem~\ref{thmmPoi2} holds for any
density in class $\cH$.
\end{cor}

%s2 ###
\section{Proofs}
\mbox{}
\begin{pf*}{Proof of Theorem \protect\ref{thmmdnstar}}
First we will show that
\begin{equation}
\limsup_{n \to\infty} d_n \leq1. \label{eqnlimsupdnstar}
\end{equation}
Fix $\ep> 0$ and let $c = 1+ \ep$. Let $n_k = k^b$, $k \geq1$, where the
constant $b$ will be chosen later. Let $W_n(c)$ be as defined in (\ref
{eqnWn}).
For each $n \geq1,$ define the events
\[
A_n := \{ W_n(c) > 0 \} .
\]
%
%$A_n$ to be the event that the graph $G(\cP_n,r_n(c, \cdot))$ has at
%least one node with out degree zero.
Recall that $\cP_n[B]$ denotes the number of points of the
point process $\cP_n$ that fall in the set $B$. Set
\[
B_k := \bigcup_{n=n_k}^{n_{k+1}} A_n,\qquad k \geq1.
\]
We will show that
\begin{equation}
\sum_{k = 1}^{\infty} P[B_k] < \infty. \label{eqnsumbk}
\end{equation}
It will then follow by the Borel--Cantelli lemma that almost surely,
only~finitely many of the events $B_k$ (and hence the $A_n$) happen.
Consequently, with probability~1, $d_n \leq1 + \ep$, eventually. Since
$\ep> 0$ is arbitrary, this will prove (\ref{eqnlimsupdnstar}).

For each $k \geq1,$ let $H_k$ be the event that there is a vertex
$X \in\cP_{n_{k+1}}$ that has out-degree zero in the graph
$G(\cP_{n_k}\cup\{X\}, r_{n_{k+1}}(c)),$ that is,
\[
H_k = \bigcup_{X \in\cP_{n_{k+1}}} \bigl\{ \cP_{n_k}[B(X,r_{n_{k+1}}(c))
\setminus\{X\}] = 0 \bigr\}.
\]
Since we have assumed the variables $N_n$ to be nondecreasing and the
functions $r_n(c,x)$ are nonincreasing
in $n$ for each fixed $c$ and $x$, we have
\[
A_n \subset H_k, \qquad n_k \leq n \leq n_{k+1}.
\]
Consequently,
\begin{equation}
B_k \subset H_k,\qquad  k \geq1. \label{eqnbk}
\end{equation}
By the Palm theory for Poisson point processes (Theorem 1.6,~\cite{Pen03}),
(\ref{eqnF}) and~(\ref{eqncutoff1}), we have
%
%e4 ###
%
\begin{eqnarray}\label{eqnhk}
P[H_k] & \leq& E\biggl[\sum_{X \in\cP_{n_{k+1}}}
1_{\{\cP_{n_k}[B(X,r_{n_{k+1}}(c))\setminus\{X\}] = 0 \}} \biggr]
\nonumber\\
& = & n_{k+1} \int_{\mR^d} e^{-n_k F(B(x,r_{n_{k+1}}(c)))} f(x) \,dx
\nonumber
\\[-8pt]
\\[-8pt]
\nonumber
& = & n_{k+1} \exp\biggl( - c \frac{n_k}{n_{k+1}} \log n_{k+1}
\biggr)
\\
& = & (k+1)^b \exp\biggl( - (1+\ep) \biggl(\frac{k}{k+1} \biggr)^b
\log(k+1)^b
\biggr).\nonumber
\end{eqnarray}
Choose $\ga> 0$ so that $(1-\ga)(1+\ep) > 1$, and pick $b$ such that
\[
\bigl((1-\ga)(1+\ep) - 1\bigr)b > 1.
\]
For sufficiently large $k$, we
have
\[
\biggl(\frac{k}{k+1} \biggr)^b > (1- \ga).
\]
Using this in (\ref{eqnhk}) , we get for all $k $ sufficiently large
\[
P[H_n] \leq\frac{1}{(k+1)^{((1-\ga)(1+\ep) - 1)b}},
\]
which is summable in $k$. (\ref{eqnsumbk}) now follows from the
above inequality and~(\ref{eqnbk}). This proves
(\ref{eqnlimsupdnstar})
by the arguments following (\ref{eqnsumbk}). To
complete the proof, we need to show that
\begin{equation}
\liminf_{n \to\infty} d_n \geq1. \label{eqnliminfdnstar}
\end{equation}
Fix $c < 1$ and pick $u$ such that $c < u < 1$. Choose $x_0$
such that $f(x_0) > 0.$ Since $f$ is continuous, we can and do fix
a $R > 0$ satisfying $g_0 u < f_0$, where
\[
f_0 = \inf_{x \in B(x_0,R)} f(x),\qquad
g_0 = \sup_{x \in B(x_0,R)} f(x).
\]
Let $\ep>0$ be such that
\begin{equation} \ep^{{1}/{d}} + c^{{1}/{d}} < u^{{1}/{d}}. \label{eqnep}
\end{equation}
Choose $R_0 > 0$ such that $2R_0 < R$ and let
$B_0 := B(x_0,R_0).$ Define the sequence of functions
$\{ \bar{r}_n(\cdot) \}_{n \geq1}$ by
\[
\bar{r}_n(v)^d = \frac{v \log n}{\th_d f_0 n},\qquad 0 \leq v \leq1.
\]
Let $\sg_n$ be the maximum number such that there exists $\sg_n$ many
disjoint balls of radius $ \bar{r}_n(u)$ with centers in
$B_0$. Then (see Lemma 2.1,~\cite{Pen99}), we can find a constant~$c_1$ such that for all $n$ sufficiently large,
\begin{equation}
\sg_n \geq\frac{c_1 n}{\log n}. \label{eqnsg}
\end{equation}
Let $\{x_1,x_2, \ldots, x_{\sg_n} \}$ be the deterministic set of points
in $B_0$ such that the balls $B(x_i,\bar{r}_n(u))$, $i=1,2, \ldots,\sg_n$,
are disjoint. Let $E_n(x)$ be the event that there is exactly one point
of $\cP_n$ in $B(x,\bar{r}_n(\ep))$ with no other point in
$B(x,\bar{r}_n(u))$, that is,
\begin{equation}\quad
E_n(x) = \{ \cP[B(x,\bar{r}_n(\ep))] = 1, \cP
_n[B(x,\bar{r}_n(u))\setminus B(x,\bar{r}_n(\ep))] = 0 \}. \label{eqnEn}\vadjust{\goodbreak}
\end{equation}
Note that the events $\{\cP[B(x,\bar{r}_n(\ep))] = 1\}$ and
$\{\cP[B(x,\bar{r}_n(u))\setminus B(x,\bar{r}_n(\ep))] = 0 \}$ are
independent. Hence for any $x \in B_0$, we have
\begin{eqnarray*}
P[E_n(x)] & = & n F(B(x,\bar{r}_n(\ep))) e^{- n F(B(x,\bar{r}_n(\ep)))}
e^{- n F(B(x,\bar{r}_n(u)))\setminus B(x,\bar{r}_n(\ep))} \\
& = & n F(B(x,\bar{r}_n(\ep))) e^{- n F(B(x,\bar{r}_n(u)))}.
\end{eqnarray*}
Note that for all $x \in B_0$ and all $n$ sufficiently large,
$B(x,\bar{r}_n(u)) \subset B(x_0,R)$. Hence for all $n$
sufficiently large and $x \in B_0$, we have
%
%
%e5 ###
%
\begin{eqnarray}\label{eqnlbEn}
P[E_n(x)] & \geq& n f_0 \th_d \bar{r}_n(\ep)^d \exp( - n g_0
\th
_d \bar{r}_n(u)^d
) \nonumber\\
& = & \ep\log n \exp\biggl( - \frac{g_0 u}{f_0} \log n \biggr)
\\
& = & \ep n^{- {g_0u}/{f_0}} \log n.\nonumber
\end{eqnarray}
Using the fact that the events $E_n(x_i)$, $i=1,\ldots,\sg_n$,
are independent, the inequality $1 - x \leq e^{-x}$, and
(\ref{eqnsg}), (\ref{eqnlbEn}), we get
%
%P[ ( \cup_{i=1}^{\sg_n} E_n(x_i) )^c ]
%& \leq& \exp( - \ep\sg_n n^{- \f{g_0u}{f_0}} \log n
%) \\
%& \leq& \exp( - c_1 \ep n^{1 - \f{g_0u}{f_0}} ),
%
\[
P\Biggl[ \Biggl( \bigcup_{i=1}^{\sg_n} E_n(x_i) \Biggr)^c \Biggr]
\leq\exp\bigl( - \ep\sg_n n^{- {(g_0u)}/{f_0}} \log n
\bigr)
\leq\exp\bigl( - c_1 \ep n^{1 - {(g_0u)}/{f_0}} \bigr),
\]
which is summable in $n$ since $g_0 u < f_0.$ Hence, by the Borel--Cantelli
lemma, almost surely, for all sufficiently large $n$ the event $E_n(x_i)$
happens for some $i = i(n)$. Hence w.p. 1, for all n sufficiently large,
we can find a random sequence $j(n)$ such that $X_{j(n)} \in\cP_n$, and
there is no other point of $\cP_n$ within a distance $\bar{r}_n(u) -
\bar{r}_n(\ep)$ of $X_{j(n)}$. By (\ref{eqnep}),
\[
\bar{r}_n(u) - \bar{r}_n(\ep) \geq\bar{r}_n(c).
\]
Since $X_{j(n)} \in B_0$, from (\ref{eqncutoff1}) and the remark above
(\ref{eqnlbEn}), we get for sufficiently large $n$,
\begin{eqnarray}\label{eqnrnbound}\qquad
\int_{B(X_{j(n)}, r_n(c,X_{j(n)}))} f(y) \,dy &=&
c \frac{\log n}{n}
\nonumber
\\[-8pt]
\\[-8pt]
\nonumber
& = &\th_d f_0 \bar{r}_n(c)^d
\leq \int_{B(X_{j(n)},\bar{r}_n(c))} f(y) \,dy.
\end{eqnarray}
Hence $r_n(c,X_{j(n)}) \leq\bar{r}_n(c)$. It follows that there
is no other point of $\cP_n$ in $B(X_{j(n)},r_n(c,X_{j(n)}))$, that is,
$X_{j(n)}$ has out-degree zero in $G_n= G(\cP_n, r_n(c)).$ Consequently
w.p. 1, $d_n \geq c$ for all $n$ sufficiently large. Since $c<1$,
this proves~(\ref{eqnliminfdnstar}).
\end{pf*}

\begin{pf*}{Proof of Theorem \protect\ref{thmmtdn}}
Since the graph $\tG_n$ is
obtained by making all the edges in $G_n$ bi-directional,
$\tW_n \leq W_n$. Consequently $\td_n \leq d_n$, and hence by
Theorem~\ref{thmmdnstar} we have
\[
\limsup_{n \to\infty} \td_n \leq1.
\]
Thus it suffices to show that
\begin{equation} \liminf_{n \to\infty} \td_n \geq1. \label{eqnliminftdn}
\end{equation}
Fix $c<1$, and let $x_0, R, R_0$ and $\bar{r}_n$ be as in the second part
of the proof of Theorem~\ref{thmmdnstar}. Recall the random
sequence $j(n)$,
defined on a set of probability one, such
that for sufficiently large $n$, the point $X_{j(n)} \in\cP_n$
has no other point of $\cP_n$ within a distance $\bar{r}_n(c)$.
Since $2R_0 < R$, by the same arguments as in (\ref{eqnrnbound}), we
have $r_n(x,c) \leq\bar{r}_n(c)$ for all $x \in B(x_0, 2 R_0)$
and all~$n$ sufficiently large.
%As in (\ref{eqnrnbound}), for any $x \in\overline{B(x_0,R)}$ (the
%closure of $B(x_0,R))$, we have $r_n(c,x) \leq
Thus almost surely, none of the points $X \in\cP_n$,
that fall in the ball $B(x_0,2R_0) \setminus B(x_0, R_0)$, have an
out-going edge to $X_{j(n)}$ for all sufficiently large $n$.

On the other hand, for any point $x \notin B(x_0,2R_0)$, suppose
$B(x,r_n(c,x)) \cap B(x_0,R_0) \ne\phi$.
%Recall that $B(x_0,2R_0) \subset B(x_0,R).$
Then we can find a point $y$ such that $\Vert y - x_0\Vert = 3R_0/2$
such that
$B(y,R_0/2) \subset B(x,r_n(c,x)).$ However,
\[
c \frac{\log n}{n} = \int_{B(x,r_n(c,x))} f(u) \,du
\geq\int_{B(y,R_0/2)} f(u) \,du \geq\frac{f_0 \th R_0^d}{2^d},
\]
which clearly is not possible for sufficiently large $n$. Thus there
can be no edge leading from a point of $\cP_n$ in $B(x_0,2R_0)^c$
to $X_{j(n)}.$ Hence the points $X_{j(n)} \in\cP_n$ have zero
in-degree as well
for sufficiently large $n$.
Consequently, with probability~1, $\td_n \geq c$ for all sufficiently
large $n$ for any
$c < 1.$ This proves (\ref{eqnliminftdn}) and thus completes the
proof of Theorem~\ref{thmmtdn}.
\end{pf*}

\begin{pf*}{Proof of Theorem \protect\ref{thmmconnectivity}}
For any $\ep> 0$, let
$c_n=c_n(\ep, \cdot)$, $n \geq1$, be as defined in (\ref{eqncn}).
Consider the random geometric graph $\bar{G}_n$ induced by the
mapping~$h$ and the enhanced graph $\tG_n = \tG_n(f,r_n(c_n,\cdot))$.
The vertex set of the graph
$\bar{G}_n$ is the set $\{h(X)\dvtx X \in\cP_n\}$ with edges between
any two
vertices $Y_i=h(X_i)$ and $Y_j = h(X_j)$ provided there is an edge between
$X_i$ and $X_j$ in the graph $\tG_n$. The vertices of $\bar{G}_n$ are
distributed according to a homogenous Poisson point process on
$[0,1]^d$ with
intensity $n$.

Let $m$ be as defined in (\ref{eqnm}).
Now by Theorem 13.2,~\cite{Pen03} (with $k_n \equiv0$), we have
%and hence $b = a_j = 0$,
%$0 \leq j \leq d-1$.
%Since $H(0) = 1$, (13.4)~\cite{Pen03}, yields
%
\[
\lim_{n \to\infty} \frac{n \th T_n^d}{\log n} = m,
\]
almost surely, where $T_n$ is the threshold for simple connectivity
in the usual uniform random geometric graph on $[0,1]^d$ (nodes being
distributed according to a homogenous Poisson point process with
intensity $n$).

By definition of $c_n(\ep)$, in the graph $\bar{G}_n$, each vertex is
connected to all its neighbors that are within a distance $m_n(\ep)$
almost surely for all sufficiently large $n$. It follows that almost surely,
the graphs $\bar{G}_n$, and hence the graphs~$\tG_n$ are connected
for all
sufficiently large $n$.
\end{pf*}

\begin{pf*}{Proof of Theorem \protect\ref{thmmmaxmindeg}}
We first show (\ref
{eqnDeln}).
Fix $c > 0$ and $\ep\in(0,1).$ Define the sequence
\[
c_n = (1+ \ep) c H_+^{-1} \biggl( \frac{1+\ep}{c} \biggr)
\log n,\qquad
n \geq1.
\]
Fix a constant $b$ such that $b \ep> 1$, and let $n_k = k^b,$ $k \geq1.$
Define the events $A_n : = \{ \Del_n \geq c_n \}$, $n \geq1$. Let
\begin{equation} B_k := \bigcup_{n=n_k}^{n_{k+1}} A_n,\qquad k \geq1.
\label{eqndefBk}
\end{equation}
Since $N_n$, $c_n$ are increasing, and $r_n(c,x)$ is decreasing
pointwise in $n$, we have
\begin{equation} B_k \subset\bigcup_{X \in\cP_{n_{k+1}}} \{ \cP
_{n_{k+1}}[B(X,r_{n_k} (c,X))] \geq c_{n_k} + 1 \}. \label{eqnmaxBk}
\end{equation}
By the Palm theory for Poisson point processes (Theorem 1.6,~\cite{Pen03}),
%
%e6 ###
%
\begin{eqnarray}\label{eqnmaxBbd}
P(B_k) & \leq& E\biggl[\sum_{X \in\cP_{n_{k+1}}}
1_{\{\cP_{n_{k+1}}[B(X,r_{n_{k}}(c,X))\setminus\{X\}]
\geq c_{n_k} \}} \biggr]
\nonumber
\\[-8pt]
\\[-8pt]
\nonumber
& = & n_{k+1} \int_{\mR^d} P \bigl( \operatorname{Po}( n_{k+1}
F(B(x,r_{n_k}(c,x))) ) \geq c_{n_k} \bigr) f(x) \,dx,
\end{eqnarray}
where $\operatorname{Po}(\lam)$ denotes a Poisson random variable with mean $\lam.$
From (\ref{eqncutoff1}), we have for sufficiently large $k$,
\[
n_{k+1} F(B(x,r_{n_k}(c,x))) = \frac{n_{k+1}}{n_k} c \log n_k
\leq(1+\ep) c b \log k.
\]
Hence for sufficiently large $k$, by definition of $H_+^{-1}$, we get
\begin{eqnarray}\label{eqnratio} \frac{c_{n_k}}{n_{k+1} F(B(x,r_{n_k}(c,x)))} &\geq&
\frac{(1+\ep)c H_+^{-1} ( {(1+\ep)}/{c} ) b \log
k} {(1+\ep) c b \log k}
\nonumber
\\[-8pt]
\\[-8pt]
\nonumber
& =& H_+^{-1}\biggl ( \frac{1+\ep}{c} \biggr)
> 1.
\end{eqnarray}
Hence using the Chernoff bound
for the Poisson distribution (see Lem\-ma~1.2,~\cite{Pen03}), we get
\begin{eqnarray*}
&&P \bigl( \operatorname{Po}( n_{k+1} F(B(x,r_{n_k}(c,x))) ) \geq c_{n_k}
\bigr)\\
&&\qquad \leq e^{ - n_{k+1} F(B(x,r_{n_k}(c,x)))
H ( {c_{n_k}}/{(n_{k+1} F(B(x,r_{n_k}(c,x))))} ) }.
\end{eqnarray*}
Since $H$ is increasing in $[1, \infty)$ and $n_{k+1} > n_k$, using
(\ref{eqnratio}) we can
bound the probability on the left-hand side in the above equation by
\[
\exp\biggl( - \frac{n_{k+1}}{n_k} c b \log k H \biggl(
H_+^{-1} \biggl( \frac{1+\ep}{c} \biggr) \biggr) \biggr)
\leq\exp\bigl( -(1+\ep) b \log k\bigr).\vadjust{\goodbreak}
\]
Substituting this bound in (\ref{eqnmaxBbd}), we get for
sufficiently large $k$,
\[
P(B_k) \leq\frac{(k+1)^b}{k^{b(1+\ep)}} \leq(1+\ep) \frac
{1}{k^{b\ep}},
\]
which is summable in $k$, since $b \ep> 1.$ Hence, by the Borel--Cantelli
lemma, almost surely only finitely many of the events $B_k$ and hence
$A_n$ happen. Hence, almost surely,
\[
\frac{\Del_n}{\log n} \leq
(1+ \ep) c H_+^{-1}\biggl ( \frac{1+\ep}{c} \biggr),
\]
eventually. The result now follows since $\ep> 0$ is arbitrary.

The result for $\delta_n(c)$ in case $c < 1$ follows from
Theorem~\ref{thmmdnstar}. For the case $c > 1,$
the proof of (\ref{eqndeln}) is entirely analogous to that of
(\ref{eqnDeln}), and so we provide the corresponding expressions.
Fix $\ep\in(0,1)$ such that $(1+ \ep) < (1- \ep) c$.
Let $b > 0$ be such that $b \ep> 1$, and define the sequence $n_k = k^b$,
$k \geq1$. Define the events $A_n := \{ \delta_n \leq c_n\}$ where
\[
c_n = (1 - \ep) c H_-^{-1} \biggl( \frac{1+\ep}{(1 - \ep)c}
\biggr)
\log n,\qquad n \geq1.
\]
Let the events $B_k$ be as defined in (\ref{eqndefBk}). The
expression analogous to
(\ref{eqnmaxBk}) will be
\begin{equation} B_k \subset\bigcup_{X \in\cP_{n_{k}}} \{ \cP
_{n_{k}}[B(X,r_{n_{k+1}} (c,X))] \leq c_{n_{k+1}} + 1 \}. \label{eqnminBk}
\end{equation}
For sufficiently large $k$,
\[
n_{k} F(B(x,r_{n_{k+1}}(c,x))) = \frac{n_{k}}{n_{k+1}} c \log n_{k+1}
\geq(1-\ep) c b \log(k+1).
\]
Hence for sufficiently large $k$, by our choice of $\ep$, we get
\begin{eqnarray} \label{eqnratio1}
\frac{c_{n_{k+1}}}{n_{k} F(B(x,r_{n_{k+1}}(c,x)))}
&\leq&\frac{(1-\ep)c H_-^{-1} ( {(1+\ep)}/{((1 - \ep)c)}
) b \log(k+1)} {(1-\ep) c b \log(k+1)}
\nonumber\hspace*{-35pt}
\\[-4pt]
\\[-12pt]
\nonumber
& =& H_-^{-1} \biggl(
\frac{1+\ep}{(1 - \ep)c} \biggr) < 1.\hspace*{-35pt}
\end{eqnarray}
Again using the Chernoff bound and proceeding as in the previous proof, we
will get
\begin{eqnarray*}
P(B_k) & \leq& n_k \exp\biggl( - (1-\ep) c b \log(k+1) \frac
{1+\ep}{(1
- \ep)c}
\biggr), \\
& \leq& (1+ \ep) \frac{1}{(k+1)^{b \ep}},
\end{eqnarray*}
which is summable in $k$. Since $\ep> 0$ is arbitrary,
(\ref{eqndeln}) now follows by the Borel--Cantelli lemma
and the arguments used earlier to infer (\ref{eqnDeln}).\vadjust{\goodbreak}
\end{pf*}

The proof of Theorem~\ref{thmmPoicgs} uses the following
lemma, which is a straightforward extension of Theorem 6.7,~\cite{Pen03}.
Let $d_{\mathrm{TV}}$ denote the total variation distance between two
random variables.
Let $\hr_n(\be), \hG_n$ and $\hW_n, $ be as
in Theorem~\ref{thmmPoicgs}. We will use
the notation, $B_n(x) = B(x,\hr_n(x))$,
\begin{eqnarray*}
\bB_n(x) &=& \bigl\{ y\dvtx \Vert y-x\Vert \leq
3 \max\{ \hr_n(x), \hr_n(y) \} \bigr\},
\\
\hat{B}_n(x) &=& \bigl\{ y\dvtx
\max\{ \hr_n(x), \hr_n(y) \} \leq\Vert y-x\Vert \leq
3 \max\{ \hr_n(x), \hr_n(y) \} \bigr\},
\end{eqnarray*}
and $\cP_n^x = \cP_n \cup\{ x \}.$ Recall that $\cP_n(B)$ denotes the
number of points of
$\cP_n$ that fall in the set $B$.
\begin{lem}\label{lemPoiapprox} Let $f$ be a continuous density on $\mR^d$. Then, for
the graph~$\hG_n$, we have
\begin{equation} d_{\mathrm{TV}}( \hW_n, \operatorname{Po}(E[\hW_n]) ) \leq
\min\biggl( 3, \frac{1}{E[\hW_n]} \biggr) \bigl(I_n^{(1)} + I_n^{(2)}\bigr),
\label{eqnpoiapprox}
\end{equation}
where
%
%e8 ###
%e7 ###
%
\begin{eqnarray}
\label{eqnI1} I_n^{(1)} & = & n^2 \int_{\mR^d} f(x) \,dx \int_{\bB_n(x) } f(y)
\,dy
\nonumber
\\[-8pt]
\\[-8pt]
\nonumber
&&{}\times P[ \cP_n(B_n(x)) = 0] P[ \cP_n(B_n(y)) = 0], \\
\label{eqnI2}I_n^{(2)} & = & n^2 \int_{\mR^d} f(x) \,dx \int_{\hat{B}_n(x)}
f(y) \,dy
\nonumber
\\[-8pt]
\\[-8pt]
\nonumber
&&{}\times
P[ \cP_n^y(B_n(x)) = 0 , \cP_n^x(B_n(y)) = 0].
\end{eqnarray}
\end{lem}

\begin{pf}
The proof is identical to
the proof of Theorem 6.7,~\cite{Pen03} with the following obvious
change. In the definition of dependency neighborhood, the parameter~$r$
is replaced by
the function $\operatorname{sup}_{x \in H_{mi} \cup H_{mj}} \hr_n(x)$.
\end{pf}

\begin{pf*}{Proof of Theorem \protect\ref{thmmPoicgs}}
By the Palm theory for
Poisson point processes,
\[
E[\hW_n] = n \int_{\mR^d} f(x) \,dx\, e^{- n \int_{B_n(x)} f(y)
\,dy} ,
\]
where $B_n(x)$ is as defined above Lemma~\ref{lemPoiapprox}. Using
(\ref{eqnrPoi}), we get
\[
E[\hW_n] = e^{-\be}.
\]
Hence, by Lemma~\ref{lemPoiapprox}, it suffices to show that $I_n^{(1)},
I_n^{(2)}$ converge to zero as \mbox{$n \to\infty$}. Again, using the Palm theory
and (\ref{eqnrPoi}), we get
%
%e9 ###
%
\begin{eqnarray}\label{eqnIn1}
I_n^{(1)} & = & n^2 \int_{\mR^d} f(x) \,dx \int_{\bB_n(x) } f(y)
\,dy\,e^{- n \int_{B_n(x)} f(u) \,du } e^{- n \int_{B_n(y)} f(v) \,dv }
\nonumber\hspace*{-35pt}
\\[-8pt]
\\[-8pt]
\nonumber
& = & e^{-2 \be} \int_{\mR^d} f(x) \,dx \int_{\bB_n(x)} f(y)
\,dy
\to 0\hspace*{-35pt}
\end{eqnarray}
as $n \to\infty$, by the dominated convergence theorem,
since $\hr_n(\be,x) \to0$ for each $x \in\mR^d$. Next
we will show that $I_n^{(2)} \to0$, as $n \to\infty$.
\begin{eqnarray}\label{eqnIn2} I_n^{(2)}& =& n^2 \int_{\mR^d} f(x) \,dx \int
_{\hat{B}_n(x) } f(y) \,dy\, e^{ - n \int_{B_n(x) \cup B_n(y)}
f(u) \,du }
\nonumber
\\[-8pt]
\\[-8pt]
\nonumber
&= &I_n^{(21)} + I_n^{(22)},
\end{eqnarray}
where
\begin{eqnarray*}
I_n^{(21)} & : = & n^2 \!\int_{\mR^d} f(x) \,dx\!
\int_{\hat{B}_n(x) \cap\{ y\dvtx \Vert y-x\Vert \geq\hr_n(x) + \hr
_n(y) \} }
f(y) \,dy\, e^{ - n \int_{B_n(x) \cup B_n(y)} f(u) \,du }, \\
I_n^{(22)} & : = & n^2 \!\int_{\mR^d} f(x) \,dx\!
\int_{\hat{A}_n(x) }
f(y) \,dy\, e^{ - n \int_{B_n(x) \cup B_n(y)} f(u) \,du },
\end{eqnarray*}
where $\hat{A}_n(x)$ is as defined in (\ref{eqnhAn}).
Consider the inner integral in $I_n^{(21)}$. On the
set $\hat{B}_n(x) \cap\{y\dvtx \hr_n(x) + \hr_n(y) \leq\Vert x-y\Vert
\}$, we have
$B_n(x) \cap B_n(y) = \phi$, and hence by~(\ref{eqnrPoi})
\[
\int_{B_n(x) \cup B_n(y)} f(u) \,du = \int_{B_n(x)} f(u) \,du +
\int_{B_n(y)} f(u) \,du = 2 \frac{\log n + \be}{n} .
\]
Thus, $I_n^{(21)}$ converges to zero using
the same arguments as in (\ref{eqnIn1}).
It remains to show that $I_n^{(22)} \to0$ as $n \to\infty$.
%
%||x-y|| < r_n(x) + r_n(y)\} } f(y) \,dy
% e^{ - n \int_{B_n(x) \cup B_n(y)} f(u) \,du } \to0,
%
%as $n \to\infty$.
Since $B_n(x) \cup B_n(y) = B_n(x) \cup K_n(x,y)$,
where $K_n(x,y) = B_n(y) \setminus B_n(x)$, we get from (\ref
{eqnrPoi}), (\ref{eqnc1}), that $I_n^{(22)}$ is bounded by a~constant times
\begin{equation} n^{1-\al} \int_{\mR^d} f(x) F(A_n(x)) \,dx,
\label{eqnfinalbound}
\end{equation}
where $A_n(x)$ is as defined in (\ref{eqnAn}). The expression in
(\ref{eqnfinalbound}) converges to zero as $n \to\infty$
by (\ref{eqnc2}).
\end{pf*}

\begin{pf*}{Proof of Theorem \protect\ref{thmmPoi2}}
Note that in the proof
of Theorem~\ref{thmmPoicgs}, (\ref{eqnc2}) is used only to prove
that (\ref{eqnfinalbound}) converges to
zero. Hence it suffices
to show that~(\ref{eqnfinalbound}) converges to
zero under (\ref{eqnc3}). Let $A_n(x)$ be as defined in (\ref
{eqnAn}). Then
\begin{eqnarray*}
A_n(x) & = & \bigl( A_n(x) \cap\{ \hr_n(y) \leq\hr_n(x) \}
\bigr)
\cup
\bigl(A_n(x) \cap\{ \hr_n(x) \leq\hr_n(y) \} \bigr) \\
& \subset& B(x,2\hr_n(x)) \cup\bigl(A_n(x) \cap\{ \hr_n(x) \leq
\hr
_n(y) \} \bigr).
\end{eqnarray*}
Using this, the expression in (\ref{eqnfinalbound}) is bounded by
\[
n^{1-\al} \biggl( \int_{\mR^d} f(x) F(B(x,2\hr_n(x))) \,dx
+ \int_{\mR^d} f(x) \,dx \int_{A_n(x) \cap\{ \hr_n(x) \leq
\hr
_n(y) \}}
f(y) \,dy \biggr).
\]
Applying Fubini's theorem to the second term above, we see that
the above expression is bounded by
\[
2 n^{1-\al} \int_{\mR^d} f(x) F(B(x,2\hr_n(x))) \,dx,
\]
which converges to zero by (\ref{eqnc3}).
\end{pf*}

\begin{pf*}{Proof of Corollary \protect\ref{corPoi}}
Let $B = B(0,1)$. Consider the
following two classes of densities:
\[
\cC_+ := \biggl\{ f\dvtx \mR^d \to[0,\infty), f \mbox{ continuous },
\int_{\mR^d} f(y) \,dy = 1,
\inf_{x \in \operatorname{Supp}(f)} f(x) > 0 \biggr\},
\]
where $\operatorname{Supp}(f)$ is the support of $f$.
%and $A^0$ denotes the interior of the set $A$,
Denote by $\cC_E$ the set of functions $f\dvtx B \to[0, \infty)$ such that
for some $r \in(0,1)$ and $p \in\mN$, $f(x) = 0, x \in B(0,r),$ and
$f(x) = A (\Vert x\Vert - r)^p, $ $x \in B\setminus B(0,r),$
with $\int_{B} f(y) \,dy = 1$. Let $\cC_I$ be the set of functions
$f\dvtx B \to[0, \infty)$ such that for some $p \in\mN$,
$f(x) = A (1-\Vert x\Vert)^p, $ $x \in B,$ with $\int_{B} f(y) \,dy = 1$.
%
%f(x) = 0, x \in B(0,r),
%f(x) = A ||x||^p, p \in\mN, x \in B\setminus B(0,r),
%
We first prove the result for densities $f \in\cC_+ \cup\cC_E \cup
\cC_I$.
To do this we need to verify conditions (\ref{eqnc1}) and (\ref{eqnc3}).

\textit{Step} 1. Fix $f \in\cC_+$. Let $S = \operatorname{Supp}(f)$ and
$f_* = \inf_{x \in S} f(x) > 0$
and $f^* = \sup_{x \in S} f(x)$. By (\ref{eqnrPoi}), for any $x \in S$
and for all $n$ large enough, we get
\[
f_* \th_d \hr_n(x)^d \leq\frac{\log n + \be}{n} \leq
f^* \th_d \hr_n(x)^d,
\]
or
\begin{equation}
\frac{1}{f^* \th_d}\frac{\log n + \be}{n} \leq
\hr_n(x)^d \leq\frac{1}{f_* \th_d}\frac{\log n + \be}{n}.
\label{boundsrn}
\end{equation}
To prove (\ref{eqnc1}), note that for any $y \in\hat{A}_n(x)$, we
can inscribe a ball of radius $\hr_n(y)/2$ inside $K_n(x,y)$. From this
and (\ref{boundsrn}), for all $n$ sufficiently large, we get
\begin{eqnarray*}
F(K_n(x,y)) & \geq& f_* \th_d \biggl( \frac{\hr_n(y)}{2}
\biggr)^d \\
& \geq& \biggl( \frac{f_*}{2^d f^*} \biggr) \biggl( \frac{\log n
+ \be}{n}
\biggr).
\end{eqnarray*}
This proves (\ref{eqnc1}) with $\al= \frac{f_*}{2^d f^*} < 1$.
By (\ref{boundsrn}) we have
\begin{equation} F(B(x, 2 \hr_n(x))) \leq f^* 2^d \th_d \hr
_n(x)^d = o(n^{1 - \al})\qquad \mbox{as } n \to\infty.
\end{equation}
This proves (\ref{eqnc3}). Thus the Poisson convergence result
holds for any $f \in\cC_+$.

\textit{Step} 2. Next we prove the result for $f \in\cC_E$. Proof for
$f \in\cC_I$ is similar and so we omit it. Let $B_r = B \setminus
B(0,r)$. Recall that $f(x) = 0$ over $B(0,r)$ and is of the form
$f(x) = A(\Vert x\Vert - r)^p$ over $B_r$.
%For any $\ga\in(0,1)$, on the set $B \setminus B(0, \ga)$,
%$f$ is strictly bounded
%away from zero, and hence (\ref{eqnc1}), (\ref{eqnc3}) hold by the
%same
%arguments as in Step 1. This is done
%to avoid complications with edge effects in the integrals that
%appear below and also to obtain certain bounds that will be used later.
%Thus the result will follow for any $f \in\cC_P$ by Theorem
%once we verify conditions (\ref{eqnc1}), (\ref{eqnc3}) on $B(0,

For any $x,y \in B_r$, we have by (\ref{eqnrPoi}),
\begin{equation}
F(B(y, \hr_n(y))) = F(B(x, \hr_n(x))) = \frac{\log
n + \be}{n}. \label{eqnequalityoverball}
\end{equation}
Note that the density is radially increasing, that is, $f(x) \leq f(y)$
if $\Vert x\Vert \leq\Vert y\Vert$. If $y \in\hat{A}_n(x)$, then $y
\notin B(x,\hr_n(x))$.
If $\Vert x\Vert \leq\Vert y\Vert$, then using (\ref{eqnequalityoverball})
and the monotonicity of $f$ we get
\[
F(K_n(x,y)) = F\bigl(B(y, \hr_n(y)) \setminus B(x, \hr_n(x))\bigr) \geq
\frac{1}{2}\biggl( \frac{\log n + \be}{n} \biggr).
\]
On the other hand if $y \in\hat{A}_n(x)$ and $\Vert y\Vert \leq\Vert
x\Vert$,
then by (\ref{eqnequalityoverball}) and the monotonicity of $f$ we have
$F(B(y, \hr_n(y)) \cap B(x, \hr_n(x))) \leq\frac{1}{2}(
\frac{\log n
+ \be
}{n} ).$ Hence
\[
F(K_n(x,y)) \geq\frac{1}{2}\biggl( \frac{\log n + \be}{n} \biggr).
\]
Thus (\ref{eqnc1}) holds with $\al= \frac{1}{2}$.
Next we verify (\ref{eqnc3}) over $B_r$. By (\ref{eqnrPoi}), we have
\begin{equation} \frac{\log n + \be}{n} = F(B(x,\hr_n(x))) = \int
_{B(x,\hr_n(x))\cap B_r} A (\Vert y\Vert-r)^{p} \,dy. \label{eqnFr}
\end{equation}
By changing to polar coordinates, the integral on the right-hand side
of the above equation has the bounds
\begin{equation} c_1 \hr_n(x)^{p+d} \leq\int_{B(x,\hr_n(x))\cap
B_r} A (\Vert y\Vert-r)^{p} \,dy \leq c_2 \hr_n(x)^{p+d}, \label{eqnboundr}
\end{equation}
for some positive constants $c_1,c_2.$
%The upper bound is obtained by
%integrating over the whole of $B(x,\hr_n(x))$ while the lower bound
%is obtained by considering the case when $\Vert x\Vert= r.$
From (\ref{eqnFr}) and (\ref{eqnboundr}), we get
\[
\hr_n(x)^{p+d} \leq c_1^{-1} \frac{\log n + \be}{n},
\]
and hence for some constant $c$,
\[
F(B(x, 2\hr_n(x))) \leq c_2 (2\hr_n(x))^{p+d} \leq c \frac{\log n +
\be}{n}
= o(n^{-{1}/{2}}).
\]
This proves (\ref{eqnc3}). The result now follows for any $f \in\cC
_E$ by Theorem~\ref{thmmPoi2}.

\textit{Step} 3. Let $f \in\cH$. Set $S_{0} = S \setminus S_1$ where
$S_1 = \bigcup_{i=1}^k (B(x_i,\delta_i)\setminus B(x_i,\eta_i))$.
The conditions of Theorem~\ref{thmmPoi2} hold for $S_0$ by
Step 1, and over $S_1$ by Step~2. This completes
the proof of Corollary~\ref{corPoi}.
\end{pf*}

% imsref loaded by akundreckaite, 2012-04-13 14:11:29
%

%suskaldyti doi

\printaddresses

\end{document}